\documentclass[12pt]{amsart}
\usepackage{amsthm}
\usepackage{amsmath, mathrsfs}
\usepackage{amssymb}
\usepackage{hyperref}
\usepackage{enumerate}
\usepackage{latexsym,array}
\usepackage{amsfonts}
\usepackage{shadow}

\newcommand{\bigA}{\mathcal{A}}

\newcommand{\bigU}{\mathcal{U}}

\newtheorem{Pa}{Paper}[section]
\newtheorem{Tm}[Pa]{{\bf Theorem}}
\newtheorem{La}[Pa]{{\bf Lemma}}
\newtheorem{Dn}[Pa]{{\bf Definition}}
\newtheorem{Cy}[Pa]{{\bf Corollary}}

\newtheorem{Rk}[Pa]{{\bf Remark}}
\newtheorem{Pn}[Pa]{{\bf Proposition}}

\date{}

\thanks{D. Alpay thanks the
Earl Katz family for endowing the chair which supported his
research, and the Binational Science Foundation Grant number
2010117.
}
\author{Daniel Alpay}
\address{(DA) Department of Mathematics
\newline
Ben Gurion University of the Negev \newline P.O.B. 653,
\newline
Be'er Sheva 84105, \newline ISRAEL}
\email{dany@math.bgu.ac.il}
\author{Guy Salomon}
\address{(GS) Department of Mathematics
\newline
Ben Gurion University of the Negev \newline P.O.B. 653,
\newline
Be'er Sheva 84105, \newline ISRAEL} \email{guysal@math.bgu.ac.il}
\title[On algebras which are inductive limits of Banach spaces]
{On algebras which are inductive limits of Banach spaces}
\begin{document}
\maketitle

\tableofcontents
%%%%%%%%%%%%%%%%%%%%%%%%%%%%%%%%%%%%%%%%%%%%%%%%%%%%%%%%%%%%%%%%%%%%%%%%%%%%%%%%%
%%%%%%%%%%%%%%%%%%%%%%%%%%%%%%%%%%%%%%%%%%%%%%%%%%%%%%%%%%%%%%%%%%%%%%%%%%%%%%%%%

\begin{abstract}
We introduce algebras which are inductive limits
of Banach spaces and carry inequalities which are counterparts of
the inequality for the norm in a Banach algebra. We then define
an associated Wiener algebra, and prove the corresponding version of the
well-known Wiener theorem. Finally, we consider factorization
theory in these algebra, and in particular, in the associated
Wiener algebra.
\end{abstract}

\section{Introduction}
The purpose of this paper is to establish a framework for
algebras which are inductive limits of Banach spaces and carry
inequalities which are counterparts of the inequality satisfied by
the norm in a Banach-algebra. More precisely, let $\bigA$ be an
algebra which is also the inductive limit of a family of Banach
spaces $\{X_\alpha:\alpha \in A\}$ directed under inclusion. We
call $\bigA$ a {\em strong algebra} if for any $\alpha \in A$
there exists $h(\alpha) \in A$ such that for any $\beta \geq
h(\alpha)$ there is a positive constant $A_{\beta,\alpha}$ for
which
\begin{equation} \label{SAdef}
\|ab\|_{\beta} \leq A_{\beta,\alpha}\|a\|_{\alpha}\|b\|_{\beta},
\quad \text{ and } \quad
\|ba\|_{\beta} \leq A_{\beta,\alpha}\|a\|_{\alpha}\|b\|_{\beta}.
\end{equation}
for every $a \in X_\alpha$ and $b \in X_\beta$. The case of a
Banach algebra
corresponds to the case where the set of indices $A$ is a singleton.\\

The strong algebras are topological algebras in the sense of
Naimark, i.e. they are locally convex, and the multiplication is
separately continuous (this follows from the universal property of
inductive limits; see Proposition \ref{pwcont}). If furthermore,
any bounded set in a strong algebra is bounded in some of the
$X_\alpha$,
%(e.g. if it is a dual of reflexive Fr\'echet space; for more
%examples see Corollary \ref{boundedSA}),
then the multiplication
is jointly continuous. The inequalities
\eqref{SAdef} express the fact that for any $\alpha \in A$,
each of the spaces $\{X_\beta:\beta \geq h(\alpha)\}$ ``absorbs"
$X_\alpha$ from both sides. Due to this property, one may
evaluate (with elements of $\bigA$) power series, and therefore,
consider
invertible elements; see for example Proposition \ref{invert} and
Theorem \ref{GL}.\\

In \cite{MR3029153} (see also \cite{vage1}), we studied a special family of such
algebras, which are inductive limits of ${\mathbf L}^2$ spaces of
measurable functions over a locally compact group. Examples
include the algebra of germs of holomorphic functions at the
origin, the Kondratiev space of Gaussian stochastic distributions
(see \cite{MR1408433} for the latter), the algebra of functions
$f:[1,\infty) \to \mathbb C$ for which $f(x)/x^p$ belongs to
${\mathbf L }^2([1,\infty))$ (this gives relations with the
theory of Dirichlet series), and a new space of non-commutative
stochastic distributions; see
\cite{MR3038506} for the latter.\\

In this paper, we first develop the general theory of strong
algebras. Then, we associate to every strong algebra a Wiener
algebra of functions in the following sense: Let $\bigA =
\underrightarrow{\lim} X_\alpha$ be a strong algebra, and define
$Y_\alpha$ to be the space of periodic functions
\[
a(t)=\sum_{n\in\mathbb Z} a_ne^{int},\quad a_n \in X_\alpha,
\]
on $-\pi \leq t <\pi$ to $\bigA$, with
\[
\|a\|_\alpha=\sum_{n\in\mathbb Z}\|a_n\|_\alpha<\infty.
\]
We call the inductive limit $\underrightarrow{\lim}Y_\alpha$ of
the Banach spaces $Y_\alpha$, {\em the Wiener algebra associated
to $\bigA$}. This family extends the case of Wiener algebras of
functions with values in a Banach
algebra. See \cite{ggk2} for the latter.\\

After showing that a Wiener algebra associated to a strong algebra is a strong
algebra itself, we prove a strong algebra counterpart of the well known theorem of
Wiener, namely, we show that an element is left/right/two-sided invertible in the
Wiener algebra associated to $\bigA$, if and only if its evaluation in any point of
the circle $t \in [-\pi,\pi)$ is left/right/two-sided invertible in $\bigA$.\\

Finally, we consider factorization theory in strong algebras (and
in particular in the associated Wiener algebra). More precisely,
we show that if $\bigA=\bigA^+ \oplus \bigA^-$, where
$\bigA^+,\bigA^-$ are closed subalgebras, then any element $a$,
which is ``close enough" (in an appropriate sense) to $1$, admits
a factorization
\[
a=a_-a_+,
\]
where $a_-,a_+$ are invertible and satisfying $a_+-1,a_+^{-1}-1 \in \bigA^+$,
$a_- -1,a_-^{-1}-1 \in \bigA^-$.\\

The paper consists of five sections besides the introduction, and
we now describe its content. A review on inductive limits,
bornological spaces, and barreled spaces, and the definition of a
strong algebra, are given in Section \ref{sec:2}. Some
topological results are given in Section \ref{sec:3}. In Section
\ref{sec:4} we study the underlying functional calculus. In
Section \ref{sec:5} we study the Wiener algebra associated to a
strong algebra, and prove the strong algebra version of the
well-known Wiener theorem. The factorization theory is given in
Section \ref{sec:6}.

%%%%%%%%%%%%%%%%%%%%%%%%%%%%%%%%%%%%%%%%%%%%%%%%%%%%%%%%%%%%%%%%%%%%%%%%
%%%%%%%%%%%%%%%%%%%%%%%%%%%%%%%%%%%%%%%%%%%%%%%%%%%%%%%%%%%%%%%%%%%%%%%%

%%%%%%%%%%%%%%%%%%%%%%%%%%%%%%%%%%%%%%%%%%%%%%%%%%%%%%%%%%%%%
%%%%%%%%%%%%%%%%%%%%%%%%%%%%%%%%%%%%%%%%%%%%%%%%%%%%%%%%%%%%%

\section{Topological review and the definition of a strong algebra}
\setcounter{equation}{0}
\label{sec:2}
To introduce strong
algebras, we first recall the definition of an inductive limit of
normed spaces. This definition can be extended to an inductive
limit of locally convex spaces; see \cite[II.27, Proposition 4; II.29, Example II]{MR83k:46003}.

\begin{Dn}
Let $\{X_\alpha:\alpha \in A\}$ be a family of subspaces of a vector space $X$ such that $X_\alpha\neq X_\beta$ for $\alpha \neq \beta$, directed under inclusion, satisfying $X=\bigcup_\alpha X_\alpha$, where
$A$ is directed under $\alpha \leq \beta$ if $X_\alpha \subseteq X_\beta$. Moreover, on each $X_\alpha$ ($\alpha \in A$),
a norm $\|\cdot\|_\alpha$ is given, such that whenever $\alpha \leq \beta$, the topology induced by $\|\cdot\|_\beta$ on $X_\alpha$ is coarser than the topology induced by $\|\cdot\|_\alpha$.
Then $X$, topologized with the inductive limit topology is
called {\em the inductive limit of the normed spaces}
$\{X_\alpha:\alpha \in A\}$.
\end{Dn}

The inductive limit has the following universal property. Given
any locally convex space $Y$, a linear map $f$ from $X$ to $Y$ is
continuous if and only if each of the restrictions
$f|_{X_\alpha}$ is continuous with respect to the topology of
$X_\alpha$; see \cite[II.27, Proposition 4]{MR83k:46003}. This property allows to take
full advantage of the inequalities \eqref{eq:ineq} in the
definition of a strong algebra given now.

\begin{Dn}
Let $\{X_\alpha:\alpha \in A\}$ be a family of Banach spaces directed under inclusions, and let $\bigA=\bigcup X_\alpha$ be its inductive limit. We call $\bigA$ a strong algebra if it is an algebra satisfying the property that for any $\alpha \in A$ there exists $h(\alpha) \in A$
such that for any $\beta \geq h(\alpha)$ there is a positive constant $A_{\beta,\alpha}$
for which
\begin{equation}
\label{eq:ineq}
\|ab\|_{\beta} \leq A_{\beta,\alpha}\|a\|_{\alpha}\|b\|_{\beta}, \quad \text{ and } \quad
\|ba\|_{\beta} \leq A_{\beta,\alpha}\|a\|_{\alpha}\|b\|_{\beta}.
\end{equation}
for every $a \in X_\alpha$ and $b \in X_\beta$.
\end{Dn}

Since a strong algebra is an inductive limit of Banach spaces, it inherits two special structures of locally
convex spaces, namely being bornological and barrelled. We recall
that a locally convex space $X$ is called bornological if every balanced, convex subset $U \subseteq X$ that
absorbs every bounded set in $X$ is a neighborhood of $0$. Equivalently, a bornological space is a locally convex
space on which each semi-norm that is bounded on bounded sets, is continuous.
We also recall that a topological vector space is said to be barrelled if each convex, balanced, closed and absorbent set is a neighborhood of zero. Equivalently, a barreled space is a locally convex space on which each semi-norm that is
semi-continuous from below, is continuous. With these definitions at hand, we can now state:

\begin{Pn}
An inductive limit of Banach spaces (and in particular a strong algebra) is bornological and barrelled.
\end{Pn}

\begin{proof}[Proof]
A Banach space is clearly barrelled and bornological, and these
two properties are kept under inductive limits;  see
\cite[III.25, Corollary1, Corollary 3; III.12 Examples 1,3 ]{MR83k:46003}.
\end{proof}

%every bornological space $X$ is the inductive limit of a
%family of normed spaces (and of Banach spaces if $X$ is quasi-complete).

%\begin{Cy}
%Every quasi-complete bornological space is barreled.
%\end{Cy}

%%%%%%%%%%%%%%%%%%%%%%%%%%%%%%%%%%%%%%%%%%%%%%%%%%%%%%%%%%%%%

\section{Topological results}
\setcounter{equation}{0}
\label{sec:3}

The term ``topological algebra" is sometimes refer to topological
vector space together with a (jointly) continuous multiplication
$(a,b) \mapsto ab$. However, in his book \cite{MR0438123}, M.A.
Naimark gives the following definition for a topological algebra.
\begin{Dn}[M.A Naimark]
$\bigA$ is called topological algebra if:
\begin{enumerate}[(a)]
\item $\bigA$ is an algebra;
\item $\bigA$ is a locally convex topological linear space;
\item the product $ab$ is a continuous function of each of the factors $a,b$ provided the other factor is fixed.
\end{enumerate}
\end{Dn}

We will show that a strong algebra is a topological algebra in the sense of Naimark.
\begin{Pn}
\label{pwcont}
Let $\bigA=\bigcup_\alpha X_\alpha$ be a strong
algebra and let $a \in \bigA$. Then the linear mappings $L_a:x
\mapsto ax$, $R_a:x \mapsto xa$ are continuous. Thus, it is
topological algebra in the sense of Naimark.
\end{Pn}
\begin{proof}[Proof]
Suppose that $a \in X_\alpha$, and
let $L_a|_{X_\gamma} : X_\beta \to \bigA$ be the restriction of the map $L_a$ to $X_\gamma$.
If $B$ is a bounded set of $X_\gamma$ then in particular we may choose $\beta \geq h(\alpha)$ such that $\beta \geq \gamma$, so
$B \subseteq \{ x \in X_\beta: \|x\|_\beta <\lambda \}$.
Thus, for any $x \in B$
\[
\|L_a|_{X_\gamma} (x) \|_\beta \leq A_{\beta,\alpha}\lambda \|a\|_\alpha.
\]
Hence, $L_a|_{X_\gamma}(B)$ is bounded in $X_\beta$ and hence in $\bigA$. Thus, for any $\gamma$, $L_a|_{X_\gamma} :X_\gamma \to \bigA$ is bounded and hence continuous, so by the universal property of the inductive limits, $L_a$ is continuous. The proof for $R_a$ is similar.
\end{proof}

The boundedness assumption in the next theorem occurs in many natural cases, which are discussed after the proof of the theorem.
\begin{Tm}\label{contSA}
If in a strong algebra $\bigA=\bigcup X_\alpha$, any set is bounded if and only if it is bounded in some of the $X_\alpha$, then the multiplication is jointly continuous.
\end{Tm}
%\begin{Cy}\label{boundedSA}
%In each of the following cases, the multiplication of the strong algebra $\bigA=\bigcup X_\alpha$ is continuous:
%\begin{enumerate}[(a)]
%\item $\bigA$ is a Banach space.
%There is $\alpha_0$ such that for any $\alpha \geq \alpha_0$ $X_\alpha=X_{\alpha_0}$ (i.e. $\bigA$ is a Banach space).
%\item $\bigA$ is a dual of a reflexive Fr\'echet space.
%\item $\bigA$ is the inductive limit of a sequence of Banach spaces $\{X_n:n \in \mathbb N\}$, and the embeddings $X_n \hookrightarrow X_{n+1}$ are compact.
%\end{enumerate}
%\end{Cy}

\begin{proof}[Proof of Theorem \ref{contSA}]
Let $(a,b) \in \bigA \times \bigA$. 
Since
\[
xy=(x-a)(y-b) + xb +ay -a b,
\]
and in view of Proposition \ref{pwcont}, $(x,y) \mapsto xy$ is
continuous in $(a,b)$ if and only if it is continuous at the
origin. So it remains to show that the multiplication is
continuous at the origin.\smallskip

Since a product of bornological spaces is bornological, $\bigA \times \bigA$ is bornological.
We will show that for every bounded set $B\subseteq \bigA \times \bigA$
and every convex, circled neighborhood $V$ of $0$ in $\bigA$, $m^{-1}(V)$ absorbs $B$
(where $m(x,y)=xy$). Hence, by the
definition of a bornological space (see Section 2) $m^{-1}(V)$ is a neighborhood of the origin.\smallskip

Thus, let $B$ be any bounded set in $\bigA \times \bigA$. Then $B
\subseteq B_1 \times B_2$ where $B_i$ $i=1,2$ are bounded sets in
$\bigA$ (see \cite[pp. 27, 5.5]{schaefer}). So there exists
$\alpha,\beta' \in A$ such that $B_1$ and $B_2$ are bounded in
$X_{\alpha}$ and $X_{\beta'}$ respectively. We may choose $\beta \geq
\beta',h(\alpha)$, and so $B$ is bounded in $X_\alpha \times
X_\beta$. In particular, $\|x\|_\alpha\|y\|_\beta$ is bounded for
all $(x,y) \in B$. Therefore, from inequality \eqref{eq:ineq},
$\|xy\|_\beta$ is bounded for all
$(x,y) \in B$, so $m(B)$ is bounded in $\bigA$. Thus, for any convex circled neighborhood $V$ of the origin, 
$V$ absorbs $m(B)$, i.e. there exists $\lambda>0$
such that $m(B) \subseteq \lambda V$. Thus $m(\sqrt{\lambda}^{-1}
B) \subseteq  V$  so $B \subseteq \sqrt{\lambda} m^{-1}(V)$.
\end{proof}

There are several natural cases in which the assumption of the previous theorem holds. Namely, in 
each of the following instances of inductive limit of Banach spaces, any bounded set of $\bigcup X_\alpha$ is 
bounded in some of the $X_\alpha$. Thus, when a strong algebra $\bigA$ is of one of these forms, 
then in particular the multiplication is jointly continuous.
\begin{enumerate}[(i)]
\item The set of indices $A$ is the singleton $\{0\}$, and hence $\bigcup X_\alpha=X_0$ is a Banach space.
\item  The set of indices $A$ is $\mathbb N$, and $\bigcup X_n$ is the {\em strict} inductive limit of the $X_n$ (and is then called an $LB$-space), that is, for any $m \geq n$ the topology of $X_n$ induced by $X_m$, is the initial topology of $X_n$. See \cite[Theorem 6.5, pp. 59]{schaefer} and \cite[III.5, Proposition 6]{MR83k:46003}.
\item The set of indices $A$ is $\mathbb N$, and the embeddings $X_n \hookrightarrow X_{n+1}$ are compact.
See \cite[III.6, Proposition 7]{MR83k:46003}.
\item The set of indices $A$ is $\mathbb N$, and the inductive limit is a dual of reflexive Fr\'echet space. More precisely, let $\Phi_1 \supseteq \Phi_2 \supseteq \dots$ be a decreasing sequence of Banach spaces, and assume that the corresponding countably normed space $\bigcap \Phi_n$ is reflexive. Then, $\bigcup \Phi_n'$, the strong dual of $\bigcap \Phi_n$ is the same as the inductive limit of the spaces $\Phi_1' \subseteq \Phi_2' \subseteq \dots$ (as a topological vector space).
See \cite[IV.23, proof of Proposition 4]{MR83k:46003} and \cite[\S 5.3, pp.45-46]{GS2_english} .
\end{enumerate}

In fact, in \cite[IV.26, Theorem 2]{MR83k:46003} of N. Bourbaki, the
following theorem is proved.
\begin{Tm}
Let $E_1$ and $E_2$ be two reflexive Fr\'echet spaces, and let $G$ a locally convex Hausdorff space. For $i=1,2$, let $F_i$ be the strong dual of $E_i$. Then every separately continuous bilinear mapping $u:F_1 \times F_2 \to G$ is continuous.
\end{Tm}
This gives another proof for the continuity of
the multiplication in case (iv).\\

There are some cases where the topology on an inductive limit
(that is, the inductive topology) is the finest topology such
that the mappings $X_\alpha \hookrightarrow \bigcup X_\alpha$ are
continuous (instead of the finest locally convex topology such
that they are continuous). One example is when $X$ is the
inductive limit of a sequence of Banach spaces $\{X_n:n \in
\mathbb N\}$, and the embeddings $X_n \hookrightarrow X_{n+1}$
are compact (see \cite[III.6, Proposition 7, Lemma 1]{MR83k:46003} and case (iii) above). In this case, we have the
following sufficient condition on mappings $X \to X$ to be
continuous.
\begin{Tm}\label{finest}
Let $X$ be the inductive limit of the family $X_\alpha$, where
its topology is the finest topology such that the
mappings $X_\alpha \hookrightarrow \bigcup X_\alpha$ are continuous. Then any map (not necessarily linear)
$f$ from an open set $W \subseteq X$ to $X$ which satisfies the
property that for any $\alpha$ there is $\beta$ such that $f(W \cap X_\alpha) \subseteq X_\beta$ and $f|_{W \cap X_\alpha}$ is continuous with respect to the topologies of $W \cap X_\alpha$ at the domain and $X_\beta$ at the range, is a continuous function $W \to X$.
\end{Tm}
\begin{proof}[Proof]
Note that in this case, $U$ is open in $X$ if and only if for any $U \cap X_\alpha$ is open in $X_\alpha$ for every $\alpha$.
Let $U$ be an open set of $X$, and let $\alpha \in A$. By the assumption,
there is $\beta$ such that
$f(W \cap X_\alpha) \subseteq X_\beta$ and $f^{-1}(U) \cap X_\alpha =f|_{W\cap X_\alpha}^{-1}(U \cap X_\beta)$ is open in $W \cap X_\alpha$.
In particular, $f^{-1}(U) \cap X_\alpha$ is open in $X_\alpha$, so $f^{-1}(U)$ is open in $X$ and hence in $W$.
\end{proof}

As a corollary to Theorem \ref{finest} we will show in the sequel that
whenever a strong algebra $\bigA$ satisfies the assumption of the theorem then the set of invertible elements is open and that $a \mapsto a^{-1}$ is continuous. See
Theorem \ref{GL}.\\

There is a ``well behaved" family of linear maps from an inductive
limit of Banach spaces into itself, which we call admissible operators.
These maps are in particular continuous, and sometimes all continuous
linear maps are of this form.

\begin{Dn}
Let $X$ be the inductive limit of the Banach spaces $X_\alpha$. Then
a linear map $T:X \to X$ which satisfies the property that for any
$\alpha$ there is $\beta$ such that $T(X_\alpha) \subseteq X_\beta$ and
$T|_{X_\alpha}$ is continuous with respect to the topologies of $X_\alpha$ for
the domain and $X_\beta$ for the range, will be called an admissible
operator of $X$. For an admissible operator $T:X \to X$ we denote by
$\|T\|^\alpha_\beta$ the norm of $T|_{X_\alpha}$ when
the range is restricted to $X_\beta$, whenever it makes sense,
and otherwise we set $\|T\|^\alpha_\beta=\infty$.
\end{Dn}

By the universal property of inductive limits, we conclude:
\begin{Pn}
Any admissible operator is continuous.
\end{Pn}
\begin{Rk}{\rm
Note that if any bounded set in $X$ is bounded in some $X_\alpha$,
then any continuous linear map is admissible.}
\end{Rk}

\begin{Rk}{\rm
Note that if $\|S\|^\alpha_\beta<\infty$ and $\|T\|^\beta_\gamma<\infty$, then
\[
\|TS\|^\alpha_\gamma \leq \|T\|^\beta_\gamma \|S\|^\alpha_\beta.
\]}
\end{Rk}
\begin{Pn}\label{invSO}
Let $T:X \to X$ be a admissible operator such that there exists $\alpha$ for which $\|T\|^\alpha_\alpha <1$ then $I-T$ is invertible, and
\[
\|(I-T)^{-1}\|^\alpha_\alpha \leq \frac{1}{1-\|T\|^\alpha_\alpha}.
\]
\end{Pn}
\begin{proof}[Proof]
Not that,
\[
\|I+T+T^2+\cdots \|^\alpha_\alpha \leq \sum_{n=0}^\infty \left(\|T\|^\alpha_\alpha\right)^n =\frac{1}{1-\|T\|^\alpha_\alpha}< \infty.
\]
Moreover, one can check that
\[
(I+T+T^2+\cdots)(I-T)=(I-T)(I+T+T^2+\cdots)=1.
\]
\end{proof}

%%%%%%%%%%%%%%%%%%%%%%%%%%%%%%%%%%%%%%%%%%%%%%%%%%%%%%%%%%%%%%%%%%%%%%%%%%%%%%%%%%%%%
%%%%%%%%%%%%%%%%%%%%%%%%%%%%%%%%%%%%%%%%%%%%%%%%%%%%%%%%%%%%%%%%%%%%%%%%%%%%%%%%%%%%%

\section{Power series and invertible elements}
\setcounter{equation}{0}
\label{sec:4}

Henceforward, we assume that $\bigA$ is a unital strong algebra.
\begin{Pn}\label{power}
Assuming $\sum_{n=0}^\infty c_n z^n$ converges in the open disk with radius $R$, then for any $a \in \bigA$ such that there exist $\alpha,\beta$ with $\beta \geq h(\alpha)$ and $A_{\beta,\alpha}\| a\|_{\alpha}<R$ it holds that
\[
\sum_{n=0}^\infty c_n a^{n} \in X_\beta \subseteq \bigA.
\]
\end{Pn}

\begin{proof}[Proof]
This follows from
\[
\sum_{n=0}^\infty |c_n|\| a^{n}\|_{\beta} \leq \sum_{n=0}^\infty |c_n| (A_{\beta,\alpha}\| a\|_{\alpha})^n \|1\|_\beta<\infty.
\]
\end{proof}

\begin{Pn}
\label{invert}
Let $a \in \bigA$ be such that there exists $\alpha,\beta$ with $\beta \geq h(\alpha)$ and $A_{\beta,\alpha}\| a\|_{\alpha}<1$ then $1-a$ is invertible
(from both sides) and it holds that
\[
\|(1-a)^{-1}\|_{\beta} \leq \frac{\|1\|_\beta}{1-A_{\beta,\alpha}\|a\|_\alpha}, \quad \|1-(1-a)^{-1}\|_{\beta} \leq \frac{A_{\beta,\alpha}\|a\|_\alpha \|1\|_\beta}{1-A_{\beta,\alpha}\|a\|_\alpha},
\]
where
\[
(1-a)^{-1}=\sum_{n=0}^\infty a^n.
\]
\end{Pn}

\begin{proof}[Proof]
Due to Proposition \ref{power} we have that
\[
\sum_{n=0}^\infty a^{n} \in X_\beta \subseteq \bigA.
\]
Moreover, clearly
\[
(1-a)\left(\sum_{n=0}^\infty a^{n}\right)=\left(\sum_{n=0}^\infty a^{n}\right)(1-a)=1,
\]
and we have that
\[
\|(1-a)^{-1}\|_\beta \leq \sum_{n=0}^\infty \| a^{n}\|_{\beta} \leq \sum_{n=0}^\infty  (A_{\beta,\alpha}\| a\|_{\alpha})^n\|1\|_\beta=\frac{\|1\|_\beta}{1-A_{\beta,\alpha}\|a\|_\alpha},
\]
and
\[
\|1-(1-a)^{-1}\|_\beta \leq \sum_{n=1}^\infty \| a^{n}\|_{\beta} \leq \sum_{n=1}^\infty  (A_{\beta,\alpha}\| a\|_{\alpha})^n\|1\|_\beta
=\frac{A_{\beta,\alpha}\|a\|_\alpha\|1\|_\beta}{1-A_{\beta,\alpha}\|a\|_\alpha}.
\]
\end{proof}

\begin{Pn}\label{left_invert}
If $a\in\bigA$ has a left inverse $a' \in X_\alpha \subseteq \bigA$ (i.e. $a'a=1$), then for any $\beta \geq h(\alpha)$ and $b \in X_\beta$ such that there exists $\gamma \geq h(\beta)$  with $A_{\gamma,\beta}A_{\beta,\alpha}\|a'\|_\alpha \|b\|_\beta<1$, it holds that $a-b$ has a left inverse $(a-b)' \in X_\gamma$,
where
\[
(a-b)'=a'\sum_{n=0}^\infty(ba')^n.
\]
and
\[
\|(a-b)'-a'\|_\gamma \leq  A_{\gamma,\alpha}\|a'\|_\alpha \frac{A_{\gamma,\beta}A_{\beta,\alpha}\|a'\|_\alpha\|b\|_\beta\|1\|_\gamma}{1-A_{\gamma,\beta}A_{\beta,\alpha}\|a'\|_\alpha\|b\|_\beta}
\]

%In particular, {\bf If $A_{p,q}\leq A_{p',q'}$ whenever $p-q \geq p'-q'$}, then the set of left invertible elements in $\bigA$ is open.
\end{Pn}
\begin{proof}[Proof]
We note that
\[
A_{\gamma,\beta}\|ba'\|_\beta \leq A_{\gamma,\beta}A_{\beta,\alpha}\|a'\|_\alpha\|b\|_\beta<1.
\]
Thus, $1-ba'$ is invertible, and
\[
a'(1-ba')^{-1}(a-b)=a'(1-ba')^{-1}(1-ba')a=1.
\]
Now, note that
\[
\begin{split}
\|(a-b)'-a'\|_\gamma
&=\|a'(1-ba')^{-1} -a'\|_\gamma\\
&\leq A_{\gamma,\alpha}\|a'\|_\alpha\|(1-ba')^{-1}-1\|_\gamma\\
&\leq A_{\gamma,\alpha}\|a'\|_\alpha \frac{A_{\gamma,\beta}\|ba'\|_\beta\|1\|_\gamma}{1-A_{\gamma,\beta}\|ba'\|_\beta}\\
&\leq A_{\gamma,\beta}\|a'\|_\alpha \frac{A_{\gamma,\beta}A_{\beta,\alpha}\|a'\|_\alpha\|b\|_\beta\|1\|_\gamma}{1-A_{\gamma,\beta}A_{\beta,\alpha}\|a'\|_\alpha\|b\|_\beta}.
\end{split}
\]

%Finally, denote
%\[
%U_a=\left\{ b \in \bigA : \text{ there is $q$ such that } \|b\|_q<\frac1{A_{q+d+1,q}A_{q,r}\|a'\|_r}  \right\}.
%\]
%Clearly $U_a$ is circled and convex. Since it absorbs every bounded set of $\bigA$, it is a neighborhood of zero.
%Now, for any $b\in U$
%\[
%A_{q+d+1,q}\|ba'\|_q \leq A_{q+d+1,q} A_{q,r}\|b\|_q \|a'\|_r<1.
%\]
%So $1-a^{-1}b$ is left invertible, and therefore $a-b=a(1-a^{-1}b)$ is left invertible too. Thus, $a+U_a$ is a set of left invertible elements.
\end{proof}

We now give a corollary to Theorem \ref{finest}.
\begin{Tm}
\label{GL}
In case where the topology on $\bigA$ is the finest topology such that the mappings $X_\alpha \hookrightarrow \bigA$ are continuous, the set of invertible elements $GL(\bigA)$ is open, and $(\cdot)^{-1}: GL(\bigA) \to GL(\bigA)$  is continuous.
\end{Tm}

\begin{proof}[Proof]
Let $a \in GL(\bigA)$ and assume that $a^{-1} \in X_\alpha$.
Let $U_a$ be the set of all $b \in \bigA$ such that there exists $\beta \geq h(\alpha)$ for which
\[
\|b\|_{\beta} < \frac{1}{A_{h(\beta),\beta}A_{\beta,\alpha}\|a^{-1}\|_\alpha}.
\]
Clearly $U_a \cap X_\beta$ is open in $X_\beta$ for any $\beta$, so $U_a$ is open. Moreover, for any $b\in U_a$
\[
A_{h(\beta),\beta}\|a^{-1}b\|_{\beta} \leq A_{h(\beta),\beta}A_{\beta,\alpha}\|a^{-1}\|_\alpha \|b\|_{\beta}<1.
\]
In view of Theorem \ref{invert}, $1-a^{-1}b$ is invertible, and therefore $a-b=a(1-a^{-1}b)$ is invertible too. Thus, $a+U_a \subseteq GL(\bigA)$, and so $GL(\bigA)$ is open.\\
Now, we note that,
\[
\begin{split}
(a+b)^{-1} -a^{-1}
&=\left(a(1+a^{-1}b)^{-1}\right)-a^{-1}\\
&=(1+a^{-1}b)^{-1}a^{-1}-a^{-1}\\
&=\left((1+a^{-1}b)^{-1}-1\right)a^{-1}.
\end{split}
\]
Therefore, for any $b \in U_a$,
\[
\begin{split}
\|(a+b)^{-1} -a^{-1}\|_{h(\beta)}
&\leq A_{h(\beta),\alpha}\|a^{-1}\|_\alpha\|(1+a^{-1}b)^{-1}-1\|_{h(\beta)}\\
&\leq A_{h(\beta),\alpha}\|a^{-1}\|_\alpha\frac{A_{h(\beta),\beta}\|a^{-1}b\|_{\beta}\|1\|_{h(\beta)}}{1-A_{h(\beta),\beta}\|a^{-1}b\|_{\beta}}\\
&\leq A_{h(\beta),\alpha}\|a^{-1}\|_\alpha\frac{A_{h(\beta),\beta}A_{\beta,\alpha}\|a^{-1}\|_{\alpha}\|b\|_{\beta}\|1\|_{h(\beta)}}{1-A_{h(\beta),\beta}A_{\beta,\alpha}\|a^{-1}\|_\alpha\|b\|_{\beta}}\\
\end{split}
\]
Thus, the function
\[
u:b \mapsto (a+b)^{-1} -a^{-1}
\]
satisfies $u(U_a \cap X_{\beta}) \subseteq
X_{h(\beta)}$, and $u|_{U_a \cap X_{\beta}}$
is continuous with respect to the topologies of
$U_a \cap X_{\beta}$ at the domain and $X_{h(\beta)}$ at
the range. So by Theorem \ref{finest} it is continuous $U_a \to \bigA$.
Since $a$ was arbitrary,  $(\cdot)^{-1}: GL(\bigA) \to GL(\bigA)$  is continuous.
\end{proof}

%%%%%%%%%%%%%%%%%%%%%%%%%%%%%%%%%%%%%%%%%%%%%%%%%%%%%%%%%%%%%
%%%%%%%%%%%%%%%%%%%%%%%%%%%%%%%%%%%%%%%%%%%%%%%%%%%%%%%%%%%%%

\section{A Wiener algebra associated to a strong algebra and a strong algebra version of the Wiener theorem}
\setcounter{equation}{0} \label{sec:5}
\begin{Dn}
Let $\bigA=\bigcup X_\alpha$ be a strong algebra.
Let $Y_\alpha$ be the space of periodic functions
\[
a(t)=\sum_{n\in\mathbb Z} a_ne^{int},\quad a_n \in X_\alpha,
\]
on $-\pi \leq t <\pi$ to $\bigA$, with
\[
\|a\|_\alpha=\sum_{n\in\mathbb Z}
\|a_n\|_\alpha<\infty.
\]
The inductive limit $\bigU=\bigcup Y_\alpha$ of the Banach spaces $Y_\alpha$ is called the Wiener algebra 
associated to $\bigA$.
\end{Dn}

\begin{Rk}{\rm
Assuming
\[
\Phi_1 \supseteq \Phi_2 \supseteq \cdots\supseteq \Phi_p \supseteq \cdots
\]
is a decreasing sequence of reflexive Banach spaces, so that $\bigcap \Phi_p$ is a reflexive Fr\'echet space, and suppose that $\bigA=\bigcup \Phi_p'$, the inductive limit of the duals, is a strong algebra (in particular its inductive limit topology coincides with its strong dual topology).
We may define
\[
\Psi_p=c_0(\mathbb Z;\Phi_p),
\]
i.e. the space of all $\Phi_p$-valued sequences $(x_n)_{n \in \mathbb Z}$ which satisfy
\[
\lim_{n\to -\infty}\|x_n\|_p=0,\quad \lim_{n\to \infty}\|x_n\|_p=0.
\]
This space is a Banach space, and its dual is
\[
\Psi_p'=\ell_1(\mathbb Z;\Phi_p'),
\]
i.e. the space of all $\Phi_p'$-valued sequences $(a_n)_{n \in \mathbb Z}$ which satisfy
\[
\sum_{n\in \mathbb Z}\|a_n\|_p <\infty.
\]
(For further reading on Banach-valued sequences spaces and their
duals we refer to \cite{MR0206189} and \cite{MR0420216}). In this
case, $\bigcup \Psi_p'$ can be identified as the Wiener algebra
associated to $\bigA$, but instead of considering only the
inductive limit topology on it, we may also consider its topology
as a strong dual of $\bigcap \Psi_p$. We do not know when these
two topologies coincide.}
\end{Rk}

\begin{Pn}
$\bigU=\bigcup_\alpha Y_\alpha$ is a strong algebra.
\end{Pn}
\begin{proof}[Proof]
For any $\alpha \in A$ and $\beta \geq h(\alpha)$ and for any $a \in Y_\alpha$ and $b \in Y_{\beta}$, it holds that
\[
\begin{split}
\|ab\|_\beta
&=\sum_{n\in\mathbb Z} \left\|\sum_{m\in\mathbb Z} a_mb_{n-m}\right\|_\beta \\
&\leq \sum_{n\in\mathbb Z} \sum_{m\in\mathbb Z} A_{\beta,\alpha}\|a_m\|_\alpha\|b_{n-m}\|_\beta\\
&\leq A_{\beta,\alpha}\|a\|_\alpha\|b\|_\beta.
\end{split}
\]
Similarly, $\|ba\|_\beta \leq A_{\beta,\alpha}\|a\|_\alpha\|b\|_\beta$.
\end{proof}

Our principal result is the following theorem:
\begin{Tm}\label{main}
For any $a \in \bigU$, $a$ is left invertible if and only if $a(t)$ is left invertible for every $t$.
\end{Tm}
To prove this theorem, we follow the strategy of the papers \cite{MR4:218g} of Bochner and Phillips and 
\cite{wiener-1932} of Wiener. The proofs are adapted to the case of strong algebras.
We begin with some lemmas.

\begin{La}\label{Paris}
If $a\in Y_\alpha \subseteq \bigU$ and if $a(0)$ has a left inverse, then there exists an element
\[
b(t)=\sum_{n\in\mathbb Z} b_ne^{int}
\]
in $\bigU$ with the following properties:
\begin{enumerate}[(i)]
\item the coefficient $b_0$ has a left inverse $b_0'$, and there exists $\beta \geq h(\alpha)$ and $\gamma \geq h(\beta)$  such that
\[
A_{\gamma,\beta}A_{\beta,\alpha}\sum_{n=1}^\infty (\|b_n\|_\alpha +\|b_{-n}\|_\alpha) \|b_{0}'\|_\beta <1.
\]
\item in some interval $t\in(-\epsilon,\epsilon)$, $b(t)=a(t)$.
\end{enumerate}
\end{La}
\begin{proof}[Proof]
We follow the proof of Wiener \cite[pp. 12-14]{wiener-1932}. On the circle $[-\pi,\pi)$, we define the functions\\
%\[
%\omega_\epsilon(t)=
%\begin{cases}
%1, &|t|<\epsilon\\
%2-\frac{|t|}{\epsilon}, &\epsilon \leq|t|<2\epsilon\\
%0, &2\epsilon \leq |t|
%\end{cases}
%\]
\setlength{\unitlength}{0.8cm}
\begin{picture}(8,4)(-11,-1)%(6,4)%(-3,-2)
\put(-10,0.5)
{$
\omega_\epsilon(t)=
\begin{cases}
1, &|t|<\epsilon\\
2-\frac{|t|}{\epsilon}, &\epsilon \leq|t|<2\epsilon\\
0, &2\epsilon \leq |t|
\end{cases}
$}
\put(-2.5,0){\vector(1,0){5}}
\put(2.7,-0.1){$t$}
\put(0,-0.5){\vector(0,1){2}}

\put(-0.5,1){\line(1,0){1}}
\put(-1,0){\line(1,2){0.5}}
\put(0.5,1){\line(1,-2){0.5}}

\put(-1,-0.1){\line(0,1){0.2}}
\put(-0.5,-0.1){\line(0,1){0.2}}
\put(0.5,-0.1){\line(0,1){0.2}}
\put(1,-0.1){\line(0,1){0.2}}

\put(-1.5,-0.5){\tiny{$-2\epsilon$}}
\put(-0.75,-0.5){\tiny{$-\epsilon$}}
\put(0.85,-0.5){\tiny{$2\epsilon$}}
\put(0.45,-0.5){\tiny{$\epsilon$}}

\put(0.2,1.4){$\omega_{\epsilon}(t)$}
\end{picture}

and
\[
b_{\epsilon}(t)=\omega_\epsilon(t)a(t)+(1-\omega_\epsilon(t))a(0)=\sum b_n(\epsilon)e^{int}.
\]
Clearly, $b_\epsilon$ satisfies the property $(ii)$ for any $\epsilon$. As for property $(i)$, by 
\cite[(2.203), p. 13]{wiener-1932},
\[
b_0(\epsilon)=a_0+\sum_{n=1}^\infty(a_{n}+a_{-n})\left(1+\frac{\cos\epsilon n-\cos 2\epsilon n}{\pi n^2 \epsilon}-\frac{3\epsilon}{2\pi}\right).
\]
Assuming $a \in Y_\alpha$, then there is $\alpha' \geq h(\alpha)$ such that $a(0)' \in X_{\alpha'}$.
Since
\[
\|b_0(\epsilon)-a(0)\|_\alpha
\leq \sum_{n=1}^\infty \left(\frac{\cos\epsilon n-\cos 2\epsilon n}{\pi n^2 \epsilon}-\frac{3\epsilon}{2\pi}\right) \|a_n+a_{-n}\|_\alpha \to 0,
\]
we may choose $\beta \geq h(\alpha')$
and $\epsilon_1$ such that for any $\epsilon \leq \epsilon_1$,
\[
A_{\beta,\alpha'}A_{\alpha',\alpha}\|b_{0}(\epsilon)-a(0)\|_\alpha \|a(0)'\|_{\alpha'}<1.
\]
Thus, by Lemma \ref{left_invert}, $b_{0}(\epsilon)$ has left inverse $b_{0}(\epsilon)' \in X_\beta$, and
\[
\|b_0(\epsilon)'-a(0)'\|_\beta \leq  A_{\beta,\alpha}\|a'(0)\|_\alpha \frac{A_{\beta,\alpha'}A_{\alpha',\alpha}\|b_0(\epsilon)-a(0)\|_\alpha \|a(0)'\|_{\alpha'}\|1\|_\beta}{1-A_{\beta,\alpha'}A_{\alpha',\alpha}\|b_0(\epsilon)-a(0)\|_\alpha\|a(0)'\|_{\alpha'}} \to 0.
\]
So we may choose $\epsilon_2 \leq \epsilon_1$, such that for any $\epsilon \leq \epsilon_2$,
\[
\|b_0(\epsilon)'\|_\beta \leq \|a(0)'\|_\beta+1.
\]
Moreover, by \cite[(2.205),(2.22) and (2.23) pp. 13-14]{wiener-1932},
\[
\sum_{n=1}^\infty (\|b_n(\epsilon)\|_\alpha+ \|b_{-n}(\epsilon)\|_\alpha) \leq \sum_{n\in\mathbb Z} \|a_n\|_\alpha A_n(\epsilon),
\]
where for sufficiently small $\epsilon$
\[
A_n(\epsilon) \leq \min\left\{ \sqrt{\epsilon}(2|n|c+\frac 9\pi),\frac {15}{ \pi} \right\}
\]
(where $c$ is some constant), so
\[
\sum_{n=1}^\infty (\|b_n(\epsilon)\|_\alpha+ \|b_{-n}(\epsilon)\|_\alpha) \leq \sum_{n\in\mathbb Z} \|a_n\|_\alpha A_n(\epsilon) \to 0.
\]
Thus, we may choose $\epsilon_3 \leq \epsilon_2$ and $\gamma \geq h(\beta)$ such that for any $\epsilon \leq \epsilon_3$
\[
A_{\gamma,\beta}A_{\beta,\alpha}\sum_{n=1}^\infty (\|b_n(\epsilon)\|_\alpha+ \|b_{-n}(\epsilon)\|_\alpha) (\|a(0)'\|_\beta+1) <1.
\]
Therefore
\[
A_{\gamma,\beta}A_{\beta,\alpha}\sum_{n=1}^\infty (\|b_n(\epsilon)\|_\alpha+ \|b_{-n}(\epsilon)\|_\alpha) \|b_{0}(\epsilon)'\|_\beta <1.
\]
Thus, we conclude that there exists $\epsilon$ small enough (i.e. $\epsilon_3$) such that
$b_\epsilon$ satisfies both properties $(i)$ and $(ii)$.
\end{proof}

\begin{La}\label{Berlin}
If $b \in \bigU$ satisfies the property $(i)$ of Lemma \ref{Paris}, then $b$ has a left inverse in $\bigU$.
\end{La}
\begin{proof}[Proof]
Setting
\[
c(t)=b(t)-b_0=\sum_{n=1}^\infty\left( b_ne^{int}+ b_{-n}e^{-int}\right),
\]
we obtain that $b=c-(-b_0)$ has a left inverse $b'$, and
\[
b'=-b_0'\sum_{n=0}^\infty(-cb_0')^n.
\]
\end{proof}
We are now ready to prove the main theorem.
\begin{proof}[Proof of Theorem \ref{main}]
If $a \in \bigU$, and if for any $t$, $a(t)$ is left invertible, then by Lemmas \ref{Paris} and \ref{Berlin},
there exists for each $t$ a function $b_{t} \in \bigU$ and $\epsilon_t$ such that $(b_{t}a)|_{(t-\epsilon_t,t+\epsilon_t)}=1|_{(t-\epsilon_t,t+\epsilon_t)}$. Since the circle is compact there exist $t_1,t_2, \cdots t_n$ such that the circle is covered by $\bigcup_{i=1}^n(t_i-\epsilon_{t_i},t_i+\epsilon_{t_i})$. We can now piece the associated functions $(b_{t_i})_{i=1}^n$ to a function $a' \in \bigU$, such that $a'a=1$.
\end{proof}

\begin{Cy}
For any $a\in\bigU$, a is invertible (from both sided) if and only if $a(t)$ is invertible (from both sided) for every $t$.
\end{Cy}

%%%%%%%%%%%%%%%%%%%%%%%%%%%%%%%%%%%%%%%%%%%%%%%%%%%%%%%%%%%%%

\section{Canonical factorization in decomposing strong algebras}
\setcounter{equation}{0}
\label{sec:6}

%\section{Canonical factorization in decomposing strong algebras}
%\setcounter{equation}{0}

A decomposing strong algebra $\bigA$ is a unital strong algebra which is a direct sum
\[
\bigA=\bigA^+ \oplus \bigA^-
\]
of two closed subalgebras. The projection of $\bigA$ onto $\bigA^+$ parallel to $\bigA^-$ will be denoted by $P$ and we set $Q=I-P$.
An element $a \in GL(\bigA)$ is said to admit a canonical factorization in case
\[
a=a_-a_+,
\]
where $a_-,a_+ \in GL(A)$ satisfy $a_+-1,a_+^{-1}-1 \in \bigA^+$, $a_- -1,a_-^{-1}-1 \in \bigA^-$.

We follow Clancey and Gohberg, who considered in \cite{MR657762}
only the case where $\bigA$ is a Banach algebra.
\begin{Tm}\label{Gohberg}
Let $\bigA$ be a decomposing strong algebra $\bigA$ in which elements that have inverse on one side are invertible. The following statements about an element $a \in \bigA$ are equivalent:
\begin{enumerate}[(a)]
\item The element $1-a$ admits a canonical factorization.
\item Each of the equations
\[
x-P(ax)=1, \quad y-Q(ya)=1
\]
is solvable in $\bigA$.
\item For any pair of elements $f,g \in \bigA$, each of the equations
\[
x-P(ax)=f, \quad y-Q(ya)=g
\]
is uniquely solvable in $\bigA$.
\end{enumerate}
\end{Tm}
This theorem and its proof is completely algebraically, so the proof given in
\cite{MR657762} still holds for the case of strong algebras.\\

Let $\bigA$ be a decomposing strong algebra. For $a \in \bigA$ we define the operators $T_a$ and $R_a$ on $\bigA$ by
\begin{equation}\label{Gohberg2}
T_a(x)=P(ax)+Q(x), \quad R_a(x)=P(x)+Q(xa).
\end{equation}
The result in Theorem \ref{Gohberg} asserts that the element $a$ admits a canonical factorization if and only if both $T_a$ and $R_a$ are invertible operators on $\bigA$.

\begin{Tm}\label{factor}
Let $\bigA$ be a decomposing strong algebra $\bigA$ in which elements that have inverse on one side are invertible, and let $a \in \bigA$. If there is $\alpha \in A$ and $\beta \geq h(\alpha)$ such that
\[
A_{\beta,\alpha}\|1-a\|_\alpha<\left(\max\{\|P\|^\beta_\beta,\|Q\|^\beta_\beta\}\right)^{-1}
\]
then $a$ admits a canonical factorization $a=a_+a_-$, and the factors $a_+,a_-$ may be chosen as $a_+=x^{-1}$ and $a_-=y^{-1}$, where $x$,$y$ are the solutions of
\[
P(ax)+Q(x)=1, \quad P(y)+Q(ya)=1,
\]
respectively.
\end{Tm}
\begin{proof}[Proof]
Let $T_a$ and $R_a$ be the operators defined on $\bigA$ as in \eqref{Gohberg2}.
Using the assumption, for any $\alpha$ such that $a \in X_\alpha$, and for any $\beta \geq h(\alpha)$, we obtain
\[
\|(I-T_a)x\|_\beta=\|P((1-a)x)\|_\beta \leq \|P\|^\beta_\beta \|(1-a)x\|_\beta \leq \|P\|^\beta_\beta A_{\beta,\alpha}\|1-a\|_\alpha \|x\|_\beta,
\]
and therefore,
\[
\|I-T_a\|^\beta_\beta \leq A_{\beta,\alpha}\|1-a\|_\alpha \|P\|^\beta_\beta <1.
\]
Similarly
\[
\|(I-R_a)x\|_\beta=\|Q(x(1-a))\|_\beta \leq \|Q\|^\beta_\beta \|x(1-a)\|_\beta \leq \|P\|^\beta_\beta A_{\beta,\alpha}\|1-a\|_\alpha \|x\|_\beta,
\]
and therefore,
\[
\|I-T_a\|^\beta_\beta \leq A_{\beta,\alpha}\|1-a\|_\alpha \|Q\|^\beta_\beta <1.
\]
Hence, by Proposition \ref{invSO}, $T_a$ and $R_a$ are
invertible, so $a$ admits a canonical factorization. By the proof
of Theorem \ref{Gohberg} (see \cite{MR657762}), it is obvious
that the factors have the stated form.
\end{proof}

One principal result is the following corollary.
\begin{Cy}
Let $\bigA=\bigcup X_\alpha$ be a strong algebra, and let
\[
\bigU= \left\{ a(t)=\sum_{n=-\infty}^\infty a_ne^{int}: \sum_{n=-\infty}^\infty \|a_n\|_\alpha<\infty \text { for some } \alpha \right\}
\]
be the associated Wiener algebra.
Denoting
\[
\bigU^+= \left\{ a\in \bigU:  a_n =0, \forall n\geq 0 \right\},
\quad \text{ and }\quad
\bigU^-_0= \left\{ a\in \bigU:  a_n =0,\forall n<0 \right\}.
\]
Then, any $a \in \bigU$, which is close enough to the identity, in the sense that
there exists $\alpha \in A$ and $\beta \geq h(\alpha)$ such that
\[
A_{\beta,\alpha}\|1-a\|_{\alpha}<1,
\]
admits a canonical factorization with respect to $\bigU^+$ and $\bigU^-_0$.
\end{Cy}
\begin{proof}[Proof]
Since the projections $P$ of $\bigU$ onto $\bigU^+$, and $Q$ of $\bigU$ onto $\bigU^-_0$, satisfy
$\|P\|^\beta_\beta=\|Q\|^\beta_\beta=1$, for all $\beta \in A$, the result follows from Theorem \ref{factor}.
\end{proof}
%%%%%%%%%%%%%%%%%%%%%%%%%%%%%%%%%%%%%%%%%%%%%%%%%%%%%%%%%%%%%%%%%%%%%%%%%%%%%%%%%
\bibliographystyle{plain}
%\bibliography{/users/faculty/math/dany/Travaux_courants/bib/all}
%\bibliography{all}
\def\cprime{$'$} \def\lfhook#1{\setbox0=\hbox{#1}{\ooalign{\hidewidth
  \lower1.5ex\hbox{'}\hidewidth\crcr\unhbox0}}} \def\cprime{$'$}
  \def\cfgrv#1{\ifmmode\setbox7\hbox{$\accent"5E#1$}\else
  \setbox7\hbox{\accent"5E#1}\penalty 10000\relax\fi\raise 1\ht7
  \hbox{\lower1.05ex\hbox to 1\wd7{\hss\accent"12\hss}}\penalty 10000
  \hskip-1\wd7\penalty 10000\box7} \def\cprime{$'$} \def\cprime{$'$}
  \def\cprime{$'$} \def\cprime{$'$}

\end{document}